\documentclass{article}
\usepackage[english]{babel}
\usepackage{lineno,titlesec}
\usepackage{color}
\usepackage{amssymb}
\usepackage{amsmath}
\usepackage{longtable}
\usepackage{hyperref,url}

\newtheorem{theorem}{Theorem}[section]

\newtheorem{example}[theorem]{Example}

\newcommand{\srg}{strongly regular graph}
\newcommand{\CC}{\mathbb{C}}
\newcommand{\RSHCD}{\mathrm{RSHCD}}
\newcommand{\PG}{\mathit{PG}}
\newcommand{\OA}{\mathit{OA}}
\newcommand{\VO}{\mathit{VO}}
\newcommand{\NO}{\mathit{NO}}
\newcommand{\NU}{\mathit{NU}}
\newcommand{\GQ}{\mathit{GQ}}

\newcommand{\codeparameters}[5]{$#1$& $#2$& $#3$& $#4$& $#5$}
\newcommand{\srgparameters}[4]{
\parbox[b]{0.95cm}{$(#1,$}\parbox[b]{0.81cm}{$#2,$}\parbox[b]{0.71cm}{$#3,$}\parbox[b]{0.71cm}{$#4)$}~
}

\newcommand{\construction}[1]{{#1}}

\titlespacing{\paragraph}{-15pt}{5pt plus 4pt minus 2pt}{0 pt plus 2pt minus 2pt}

\newcommand{\srgparms}[4]{\paragraph{\protect \srgparameters{#1}{#2}{#3}{#4}}}
\newcommand{\srgps}[4]{$#1$& $#2$& $#3$& $#4$}

\providecommand{\keywords}[1]{\textbf{\textit{Keywords---}} #1}

\begin{document}
\title{Implementing Brouwer's database of \srg s}

\author{Nathann Cohen\thanks{CNRS and Universit\'{e} Paris-Sud 11,
\url{nathann.cohen@gmail.com}}
\and
Dmitrii V. Pasechnik\thanks{Department of Computer Science, The University of Oxford, UK,
\url{dimpase@cs.ox.ac.uk}}}

\date{}
\maketitle

\begin{abstract}
  Andries Brouwer maintains a public database of existence results for \srg s on
  $n\leq 1300$ vertices. We have implemented most of the infinite families of graphs
  listed there in the open-source software Sagemath~\cite{sage}, as well as
  provided constructions of the ``sporadic'' cases, to obtain a graph for each
  set of parameters with known examples.  Besides providing a convenient way to
  verify these existence results from the actual graphs, it also extends the
  database to higher values of $n$.
\end{abstract}

\keywords{
05E30: strongly regular graphs, association schemes;
68-04: explicit machine computation and programs
}

\section{Introduction}

Many researchers in algebraic combinatorics or an adjacent field at some
point want to get their hands on a list of feasible parameters of \srg s,
and on actual examples of graphs. These graphs are studied and/or used in
hundreds of articles; recent highlights
in using \srg s include A. Bondarenko's \cite{MR3201240} and an improvement
of the latter by T. Jenrich and A.\,E. Brouwer \cite{MR3292266}.
While parameters are available from A.\,E. Brouwer's online database \cite{aebdb},
actually constructing an example can easily take a lot of time and effort.
The project described here aims at making these tasks almost trivial by providing
the necessary graph constructions, and a way to obtain a \srg\ from a tuple of parameters,
in the computer algebra system Sagemath~\cite{sage} (also known as Sage).
It is worth mentioning that a large part of the project makes use of GAP~\cite{GAP4}
and its packages, in particular L. Soicher's package GRAPE~\cite{GRAPE04}.

Exhaustive enumeration of the non-isomorphic \srg s has been performed for some
tuples $\mathcal{T}$ of parameters (see T. Spence \cite{tsdb}). However, the sheer
number of non-isomorphic examples (see e.g. D. Fon-Der-Flaass \cite{MR1924761}
or M. Muzychuk \cite{Muzychuk07}) makes it hard to expect to be able to generate
all of them, for a given $\mathcal{T}$,
in reasonable time. Thus we opted for a minimalist approach: for
each $\mathcal{T}$ we generate an example, provided that one is known.  We
note, however, that some of constructions implemented can generate many
examples with the same $\mathcal{T}$; e.g. we have implemented the construction to generate the point graph
of the generalized quadrangle $T_2^*(\mathcal{O})$ (see \cite{PT09}) from any
hyperoval $\mathcal{O}\subset PG(2,2^k)$. As well, many $\mathcal{T}$'s can be realized by
more than one implemented construction, sometimes leading to isomorphic graphs, and sometimes not.

Our desire to take on this project was motivated by the following considerations.
\begin{itemize}
\item One wants to double-check that the constructions are correct and their descriptions
are complete; indeed, a program is more trustworthy than a proof in some situations,
and coding a construction is a good test for completeness of the description provided.
\item We wanted to see that the Sage combinatorial, graph-theoretic, and group-theoretic primitives
to deal with such constructions are mature and versatile, so that coding of constructions is
relatively easy and quick.
\item One learns a lot while working on such a project, both the underlying mathematics, and
how the toolset can be improved. In particular, one might come along simplifications of
constructions, and this actually happened on couple of occasions, see Sect.~\ref{novel}.
\item As time goes by, possible gaps in constructions are harder and harder to fill in.
Reconstructing omitted proof details becomes a tricky and time-consuming task.
\end{itemize}

In particular, as far as the latter item is concerned, we seem to have uncovered
at least one substantial gap in constructions (see
Sect.~\ref{sect:incorrect}). Furthermore, a number of constructions needed
feedback from their authors or discussions with experts -- sometimes quite substantial -- to code them.

A large part of the constructions use in a nontrivial way another combinatorial
or algebraic object: block design, Hadamard matrix, two-graph, two-distance
code, finite group, etc.  In particular, at the start of the project some of
these were lacking in Sagemath, we needed to implement constructions of certain
block designs, regular symmetric Hadamard matrices with constant diagonal (where
the gap just mentioned was uncovered), skew-Hadamard matrices, and
two-graphs. As well, we created a small database of two-distance codes (see
Sect.~\ref{sec:twoweight}).

The remainder of the paper consists of a short introduction to \srg s,
pointing out particular relevant Sagemath features, and a description of our
implementations, few of them seemingly novel.
We list the constructions that we implemented, and discuss few gaps
we discovered in the literature.

\section{Strongly regular graphs and related objects}
An undirected regular degree $k$ graph $\Gamma$ on $n$ vertices (with $0<k<n-1$)
is called \emph{strongly regular} if the vertices $u$ and $v$ of any edge have
$\lambda$ common neighbours, and the vertices $u$ and $v$ of any non-edge have
$\mu$ common neighbours. One says that $\Gamma$ has \emph{parameters}
$(n,k,\lambda,\mu)$.  Note that the \emph{complement} of $\Gamma$, i.e. the
graph with the same set of vertices and edges being precisely the non-edges of
$\Gamma$, is also a strongly regular whose parameters are related by a simple
formula to these of $\Gamma$ (see e.g. A.\,E. Brouwer and W. Haemers \cite{BH12} for details).

\begin{example}\label{johnson}
Let $\Gamma$ be the graph with vertices being $k$-subsets of an $m$-set, with
$k\leq\lfloor m/2\rfloor$; two vertices are adjacent if the corresponding $k$-subsets
intersect in a $(k-1)$-subset. Such graphs are called \emph{Johnson graphs} and
denoted by $J(m,k)$ (in Sagemath, $J(m,k)$
can be constructed by calling the function {\tt graphs.JohnsonGraph(m,k)}).
Then $J(m,2)$ is a \srg, with parameters $(\binom{m}{2},2(m-2),m-2,4)$.
\end{example}

\begin{example}\label{rshcd}
RSHCD -- a $(n,\epsilon)$-Regular Symmetric Hadamard matrix $M$ with Constant
Diagonal is an $n\times n$ symmetric $\pm 1$-matrix such that: 1)
$MM^T=nI$; 2) its rows sums are all equal to $\delta\epsilon\sqrt{n}$,
where $\epsilon\in \{-1,+1\}$ and $\delta$ is the (constant) diagonal value of
$M$, usually denoted $\RSHCD^-$ and $\RSHCD^+$.
These matrices yield regular two-graphs. As well, they yield
\srg s: replacing all the entries equal to the diagonal values
by 0, and the remaining entries by 1 gives the adjacency matrix of a \srg .
\end{example}

Some sources further require that both $\Gamma$ and its complement are connected;
in terms of parameters this means $0<\mu<k$.
This excludes the trivial case of $\Gamma$ (or its complement)
being disjoint union of complete graphs of the same size. Sagemath implementation
does not impose this restriction.

A considerable number of techniques ruling out the existence of a \srg\
$\Gamma$ with given parameters $(n,k,\lambda,\mu)$ are known, e.g. based on
computing eigenvalues of the adjacency matrix $A$ of $\Gamma$. As $A$ generates
a dimension 3 commutative subalgebra of $\CC^{n\times n}$, one sees that there
are just 3 distinct eigenvalues of $A$, and they are determined by the
parameters (e.g. the largest eigenvalue is $k$).  Sagemath implements
parameter-based techniques to rule out sets of parameters from A.\,E. Brouwer and J.\,H. van Lint \cite{BvL84},
and from A.\,E. Brouwer, A.\,M. Cohen, and A. Neumaier \cite{BCN89}.

We use standard terminology for finite permutation groups, finite simple groups, and
geometries over finite fields from \cite{BCN89,BH12}.

\section{Structure and use of the implementation}

The strongly regular graphs are split into two categories: the fixed-size graphs
(see Sect.\ref{sec:small}) and the families of strongly regular graphs (see
Sect.\ref{sec:family}). The parameters $(n,k,\lambda,\mu)$ of fixed-size graphs
are hardcoded, while each family of strongly regular graphs has a helper function
which takes $(n,k,\lambda,\mu)$ as an INPUT and answers whether the graph
family is able to produce a graph with the required parameters. Some
families forward
their queries to the databases of Balanced Incomplete Block Designs,
of Orthogonal Arrays,  of Hadamard
matrices of various types, and of 2-weight codes.

With this design, it takes $\leq 3$ seconds on a modern laptop to know which
graphs on $<1300$ vertices can be produced by the implemented constructions
(i.e. as far as the online database goes).\\

In Sage, everything is made available to the user through a single function {\tt
  graphs.strongly\_regular\_graph} that produces a graph matching the provided
parameters. Note that $\mu$, the fourth parameter, can be omitted.

\noindent \\
\phantom{eeee}{\tt sage:\hspace{1mm}G=graphs.strongly\_regular\_graph(175,30,5)}\\
\phantom{eeee}{\tt sage:\hspace{1mm}G}\\
\phantom{eeee}{\tt  AS(5)*; GQ(6, 4): Graph on 175 vertices}\\

One can use the same function to learn whether a set of parameters is realizable,
if it is not, or if the existence problem is unsettled:

\noindent \\
\phantom{eeee}{\tt sage:\hspace{1mm}graphs.strongly\_regular\_graph(175,30,5,5,existence=True)}\\
\phantom{eeee}{\tt True}\\
\phantom{eeee}{\tt sage:\hspace{1mm}graphs.strongly\_regular\_graph(57,14,1,existence=True)}\\
\phantom{eeee}{\tt False}\\
\phantom{eeee}{\tt sage:\hspace{1mm}graphs.strongly\_regular\_graph(3250,57,0,1,existence=True)}\\
\phantom{eeee}{\tt Unknown}

Technical details and descriptions of many specific functions is available as a part \cite{sage_srg} of
the Sagemath manual, which can be found online \cite{sage}.

\section{Fixed-size constructions}
\label{sec:small}

\subsection{``Sporadic'' examples}
Here we did not attempt to give an exhaustive list of references for each graph,
for some of them have several papers devoted to them in one or another way.

We identify the corresponding graphs by their parameters, and provide references
and some construction details for each of them.

\srgparms{36}{14}{4}{6} Hubaut~\cite[S.9]{Hu75}. 
\construction{Subgraph of common neighbours of a triangle in Suzuki graph.}

\srgparms{50}{7}{0}{1} \cite[Sect.9.1.7 (iv)]{BH12}. The Hoffman-Singleton graph.

\srgparms{56}{10}{0}{2} \cite[Sect.9.1.7 (v)]{BH12}. The Sims-Gewirtz graph.

\srgparms{77}{16}{0}{4} \cite[Sect.9.1.7 (vi)]{BH12}. The $M_{22}$-graph.

\srgparms{100}{22}{0}{6} \cite[Sect.9.1.7 (vii)]{BH12}. The Higman-Sims graph.

\srgparms{100}{44}{18}{20} J{\o}rgensen and Klin~\cite{JK03}. 
\construction{Built as a Cayley graph.}

\srgparms{100}{45}{20}{20} \cite{JK03}. 
\construction{Built as a Cayley graph.}

\srgparms{105}{32}{4}{12} Goethals and Seidel~\cite{GS70}, Coolsaet~\cite{Co06}.

\srgparms{120}{63}{30}{36} R. Mathon, cf. \cite[Sect.6.A]{BvL84}.
\construction{ The distance-2 graph of $J(10,3)$.}

\srgparms{120}{77}{52}{44} Unique by J. Degraer  K. Coolsaet \cite{DC08}.
\construction{We first build a $2-(21,7,12)$ design, by removing
two points from the Witt design on 23 points. We then build the
intersection graph of blocks with intersection size 3.}

\srgparms{126}{25}{8}{4} R. Mathon, cf. \cite[Sect.6.A]{BvL84}.
\construction{The distance-(1 or 4) graph of $J(9,4)$.}

\srgparms{126}{50}{13}{24} Goethals, cf. \cite{BvL84}.

\srgparms{144}{39}{6}{12}
A.\,A. Ivanov, M.\,H. Klin, and I.\,A. Faradjev~\cite[Table 9]{IKF2}.
\construction{An orbital of degree 39 (among 2 such orbitals) of the group
  $PSL_3(3)$ acting on the (right) cosets of a subgroup of order 39.}

\srgparms{162}{56}{10}{24} \cite[S.12]{Hu75}.
\construction{The complement of the subgraph induced on the neighbours
of a vertex in the complement of McLaughlin graph.}

\srgparms{175}{72}{20}{36} \cite[Sect.10.B (iv)]{BvL84}.
\construction{Obtained from the line graph $\Lambda$ of
  Hoffman-Singleton Graph, by setting two vertices to be adjacent if their
  distance in $\Lambda$ is exactly 2. For more information, see
  \url{http://www.win.tue.nl/~aeb/graphs/McL.html}.}

\srgparms{176}{49}{12}{14} Brouwer~\cite{BrouwerPolarities82}.
\construction{Built from the symmetric Higman-Sims design. There exists an involution
$\sigma$ exchanging the points and blocks of the Higman-Sims design, such
that each point is mapped onto a block that contains it (i.e. $\sigma$ is a
polarity with all absolute points). The graph is then built by making two
vertices $u,v$ adjacent whenever $v\in \sigma(u)$.}

\srgparms{176}{85}{48}{34} W. Haemers, cf. \cite[Sect.10.B.(vi)]{BvL84}.
\construction{Obtained from the $(175,72,20,36)$-graph  by attaching a isolated
  vertex and doing \emph{Seidel switching} (cf. \cite[Sect.10.6.1]{BH12})
with respect to the disjoint union of 18 maximum cliques.}

\srgparms{176}{105}{68}{54} \cite[S.7]{Hu75}; (a rank 3 representation of $M_{22}$).
\construction{We first build a $2-(22,7,16)$ design, by removing
one point from the Witt design on 23 points. We then build the
intersection graph of blocks with intersection size 3.}

\srgparms{196}{91}{42}{42} Ionin and Shrikhande~\cite{IS06}.

\srgparms{210}{99}{48}{45} Klin {\em et al.}~\cite{KPRWZ10}.
\construction{$S_7$ acts on the 210 digraphs isomorphic to the
disjoint union of $K_1$ and the circulant 6-vertex digraph which one
can obtain using Sagemath as {\tt digraphs.Circulant(6,[1,4])}.
This action has 16 orbitals; the package \cite{COCO} found a merging of them,
explicitly described in \cite{KPRWZ10}, resulting in this graph.}

\srgparms{231}{30}{9}{3} Brouwer~\cite{Br86}. The Cameron graph.

\srgparms{243}{110}{37}{60} Goethals and Seidel~\cite{GS75}.
\construction{Consider the orthogonal complement of
the ternary Golay code, which has 243 words. On them we define a
graph, with two words adjacent if their Hamming distance is 9.}

\srgparms{253}{140}{87}{65} \cite[S.6]{Hu75};  a rank 3 representation of $M_{23}$.
\construction{We first build the Witt design on 23 points which is a
$2-(23,7,21)$ design. We then build the intersection graph of blocks with
intersection size 3.}

\srgparms{275}{112}{30}{56}  \cite[S.13]{Hu75}. The McLaughlin graph.

\srgparms{276}{140}{58}{84} Haemers and Tonchev~\cite{HT96}.
\construction{The graph is built from from McLaughlin graph, with an added
isolated vertex. We then perform Seidel switching on a set of
28 disjoint 5-cliques.}

\srgparms{280}{117}{44}{52} Mathon and Rosa~\cite{MR85}.
\construction{The vertices of the graph are all 280 partitions of a set of cardinality 9 into
3-sets, e.g. $\{\{a,b,c\},\{d,e,f\},\{g,h,i\}\}$. The cross-intersection of two
partitions $P=\{P_1,P_2,P_3\}$ and $P'=\{P'_1,P'_2,P'_3\}$ being defined as
$\{P_i \cap P'_j: 1\leq i,j\leq 3\}$, two vertices of `G` are set to be adjacent
if the cross-intersection of their respective partitions does not contain
exactly 7 nonempty sets.}

\srgparms{280}{135}{70}{60} \cite[Table 9, p.51]{IKF2}.
\construction{This graph is built from the rank 4 action of $J_2$ on the cosets of a subgroup $3.PGL(2,9)$.}

\srgparms{324}{152}{70}{72} See Sect.~\ref{rshcd324}.
\construction{$\RSHCD^-$, see Example~\ref{rshcd}.  We build an apparently new example using the
(324, 153, 72, 72)-graph; other example may be found in \cite{HX10}.}

\srgparms{324}{153}{72}{72} See Sect.~\ref{rshcd324}.
\construction{$\RSHCD^+$, see Example~\ref{rshcd}. We build the example from  \cite{JKT01};
more examples may be found in \cite{MX06,HX10}.}

\srgparms{416}{100}{36}{20} \cite[S.14]{Hu75}; (rank 3 representation of $G_2(4)$).
\construction{This graph is isomorphic to the subgraph
of the Suzuki graph \cite[S.15]{Hu75} induced on the neighbours of a vertex.}

\srgparms{560}{208}{72}{80} \cite[Table 9, p.45]{IKF2}.
\construction{Obtained as the union of 4 orbitals (among the 13 that exist) of the group
$Sz(8)$ in its primitive action on 560 points.}

\srgparms{630}{85}{20}{10} W. Haemers \cite{Haemers81}, see also \cite[Sect.10.B.(v)]{BvL84}.
\construction{This graph is the line graph of a $pg(5,18,2)$; its point
graph is the $(175,72,20,36)$-srg from this table.  One then selects a subset of 630 maximum
cliques in the latter to form the set of lines of the $pg(5,18,2)$.}

\srgparms{765}{192}{48}{48} Ionin and Kharaghani, see Sect. \ref{sect:IK765}.

\srgparms{784}{243}{82}{72} R. Mathon, cf. \cite[Sect.6.D]{BvL84}.
\construction{This and the following two are Mathon's graphs from merging classes in the product of pseudo-cyclic
association scheme for action of $O_3(8)$ on elliptic lines in $\PG(2,8)$, studied
by H.\,D.\,L. Hollmann \cite{Ho82}.}

\srgparms{784}{270}{98}{90} R. Mathon, cf. \cite[Sect.6.D]{BvL84}.
\srgparms{784}{297}{116}{110} R. Mathon, cf. \cite[Sect.6.D]{BvL84}.

\srgparms{936}{375}{150}{150} Janko and Kharaghani \cite{JankoKharaghani02}.

\srgparms{1288}{792}{476}{504} Brouwer and van Eijl~\cite{BvE92}.
\construction{This graph is built on the words of weight 12 in the binary
Golay code. Two of them are then made adjacent if their symmetric difference
has weight 12.}

\srgparms{1782}{416}{100}{96} \cite[S.15]{Hu75}. Suzuki graph, rank 3 representation of $Suz$.

\srgparms{1800}{1029}{588}{588} Janko and Kharaghani \cite{JankoKharaghani02}.

\subsection{Two-weight codes database}
\label{sec:twoweight}

The rest of the fixed-size constructions of \srg s
in the database originate from linear $d$-dimensional two-weight codes of length $\ell$
with weights $w_1$ and $w_2$ over $\mathbb{F}_q$. We use data
shared by Eric Chen~\cite{ChenDB}, data by
Axel Kohnert~\cite{Kohnert07} shared by Alfred Wassermann, data from
I. Bouyukliev and J. Simonis \cite[Theorem 4.1]{BS03}, and from
L.\,A. Disset~\cite{Disset00}.
\begin{longtable}{rrrr|rrrrr|r}
\multicolumn{4}{c|}{Graph parameters} &
\multicolumn{5}{c|}{Code parameters}&Ref.\\
$n$& $k$& $\lambda$& $\mu$ & $q$& $\ell$& $d$& $w_1$& $w_2$\\
\hline\\
\endhead
\srgps{81}{50}{31}{30}&
\codeparameters{3}{15}{4}{9}{12}&\cite{ChenDB}\\
\srgps{243}{220}{199}{200}&
\codeparameters{3}{55}{5}{36}{45}&\cite{ChenDB}\\
\srgps{256}{153}{92}{90}&
\codeparameters{4}{34}{4}{24}{28}&\cite{ChenDB}\\
\srgps{256}{170}{114}{110}&
\codeparameters{2}{85}{8}{40}{48}&\cite{ChenDB}\\
\srgps{256}{187}{138}{132}&
\codeparameters{2}{68}{8}{32}{40}&\cite{ChenDB}\\
\srgps{512}{73}{12}{10}&
\codeparameters{2}{219}{9}{96}{112}&\cite{ChenDB}\\
\srgps{512}{219}{102}{84}&
\codeparameters{2}{73}{9}{32}{40}&\cite{ChenDB}\\
\srgps{512}{315}{202}{180}&
\codeparameters{2}{70}{9}{32}{40}&\cite{Kohnert07}\\
\srgps{625}{364}{213}{210}&
\codeparameters{5}{65}{4}{50}{55}&\cite{ChenDB}\\
\srgps{625}{416}{279}{272}&
\codeparameters{5}{52}{4}{40}{45}&\cite{ChenDB}\\
\srgps{625}{468}{353}{342}&
\codeparameters{5}{39}{4}{30}{35}&\cite{BS03}\\
\srgps{729}{336}{153}{156}&
\codeparameters{3}{168}{6}{108}{117}&\cite{Disset00}\\
\srgps{729}{420}{243}{240}&
\codeparameters{3}{154}{6}{99}{108}&\cite{ChenDB}\\
\srgps{729}{448}{277}{272}&
\codeparameters{3}{140}{6}{90}{99}&\cite{Kohnert07}\\
\srgps{729}{476}{313}{306}&
\codeparameters{3}{126}{6}{81}{90}&\cite{ChenDB}\\
\srgps{729}{532}{391}{380}&
\codeparameters{3}{98}{6}{63}{72}&\cite{ChenDB}\\
\srgps{729}{560}{433}{420}&
\codeparameters{3}{84}{6}{54}{63}&\cite{ChenDB}\\
\srgps{729}{616}{523}{506}&
\codeparameters{3}{56}{6}{36}{45}&\cite{ChenDB}\\
\srgps{1024}{363}{122}{132}&
\codeparameters{4}{121}{5}{88}{96}&\cite{Disset00}\\
\srgps{1024}{396}{148}{156}&
\codeparameters{4}{132}{5}{96}{104}&\cite{Disset00}\\
\srgps{1024}{429}{176}{182}&
\codeparameters{4}{143}{5}{104}{112}&\cite{Disset00}\\
\srgps{1024}{825}{668}{650}&
\codeparameters{2}{198}{10}{96}{112}&\cite{ChenDB}\\
\hline
\end{longtable}
Note that some of these codes are members of infinite families; this will be
explored and extended in forthcoming work.

\section{Infinite families}
\label{sec:family}

These are roughly divided into two parts: graphs related to
finite geometries over finite fields (in particular
various classical geometries), and graphs obtained
by combinatorial constructions.

\subsection{Graphs from finite geometries}
Here $q$ denotes a prime power, and $\epsilon\in\{-,+\}$.
\begin{itemize}
\item Graphs arising from projective geometry designs are discussed in Sect.~\ref{combgraphs},
along with other Steiner graphs.
\item Paley graphs. The vertices are the elements of $\mathbb{F}_q$, with $q\equiv 1\mod 4$;
two vertices are adjacent if their difference is a nonzero square in $\mathbb{F}_q$; see
\cite[9.1.2]{BH12}.
\item Polar space graphs. These include polar spaces for orthogonal and
unitary groups, see entries $O^\epsilon_{2d}(q)$,
$O_{2d+1}(q)$, and $U_d(q)$ in \cite[Table 9.9]{BH12}.
Sagemath also has an implementation of polar spaces for
symplectic groups (entry $Sp_{2d}(q)$ in [loc.cit.]), but we do not use them in the database, as they have the
same parameters as these for orthogonal groups.
\item Generalised quadrangle graphs, $\GQ(s,t)$ in  \cite[Table 9.9]{BH12}.
Apart from these appearing as polar space graphs, with $s=t=q$, $s^2=t=q$,
and $s=q^2$, $t=q^3$, we provide other examples, as follows.
\begin{itemize}
\item Unitary dual polar graphs. This gives $s=q^3$, $t=q^2$.
\item $\GQ(q-1,q+1)$-graphs for $q$ odd are constructed following
 Ahrens and Szekeres, see \cite[3.1.5]{PT09}, and for $q$ even
we provide the $T_2^*(\mathcal{O})$ construction, see \cite[3.1.3]{PT09},
from a hyperoval $\mathcal{O}$ in $\PG(2,q)$.
\item $\GQ(q+1,q-1)$ are constructed as line graphs of  $\GQ(q-1,q+1)$.
\end{itemize}
\item Affine polar graphs. These are the entry $\VO^\epsilon_{2d}(q)$ in  \cite[Table 9.9]{BH12}.
\item Graphs of non-degenerate hyperplanes of orthogonal polar spaces,
with adjacency specified by degenerate intersection;
see $\NO_{2d+1}^\epsilon(q)$ in \cite[Table 9.9]{BH12}.
These are constructions by Wilbrink, cf. \cite[Sect.7.C]{BvL84}.
The implementation in Sagemath simply takes the appropriate orbit and orbital of the orthogonal
group acting on the hyperplanes using parameters of the graph, namely
$v=q^d(q^d+\epsilon)/2$, $k=(q^d-\epsilon)(q^{d-1}+\epsilon)$.
\item Graphs of non-isotropic points of polar spaces, with adjacency specified by orthogonality.
These include a number of cases.
\begin{itemize}
\item Non-isotropic points of orthogonal polar spaces over $\mathbb{F}_2$; see
$\NO^\epsilon_{2d}(2)$ in \cite[Table 9.9]{BH12}.
\item One class of non-isotropic points of orthogonal polar spaces over $\mathbb{F}_3$; see
$\NO^\epsilon_{2d}(3)$ in \cite[Table 9.9]{BH12}.
\item One class of non-isotropic points of orthogonal polar spaces (specified by a non-degenerate
quadratic form $F$) over $\mathbb{F}_5$; see
$\NO^{\epsilon\perp}_{2d+1}(5)$ in \cite[Table 9.9]{BH12}. This is a construction by Wilbrink, cf.
\cite[Sect.7.D]{BvL84}, where the class of points $p$ is described in terms of the type of the
quadric specified by $p^\perp\cap Q$, where $Q$ is the set of isotropic points of the space,
i.e. $Q:=\{x\in PG(2d,5)\mid F(x)=0\}$, and $p^\perp :=\{x\in Q\mid F(p+x)=F(p)\}$.
The implementation in Sagemath takes $\{x\in PG(2d,5)\mid F(x)=\pm 1\}$ for $\epsilon=+$,
and the rest of non-isotropic points for $\epsilon=-$.
\item Non-isotropic points of unitary polar spaces; see $\NU{d}(q)$ in \cite[Table 9.9]{BH12}.
\end{itemize}
\item Graphs of Taylor two-graphs, see  \cite[Table 9.9]{BH12} and \cite[Sect.7E]{BvL84}. Note that we implement an
efficient construction that does not need all the triples of the corresponding two-graphs,
by first directly constructing the descendant \srg s on $q^3$ vertices, and a partition of its vertices into cliques. The latter
provides a set to perform Seidel switching on the disjoint union with $K_1$, and
obtain the \srg\ on $q^3+1$ vertices. See Sagemath documentation for {\tt graphs.TaylorTwographSRG}
for details.
\item Cossidente-Penttila hemisystems in $\PG(5,q)$, for $q$ odd prime
power \cite{CP05}, are certain partitions of points of the minus type quadric in $\PG(5,q)$ into
two parts $V$, $V'$ of equal size. The subgraph $\Gamma$
of the collinearity graph of the corresponding $\GQ(q,q^2)$ induced on $V$
has parameters $((q^3+1)(q+1)/2,(q^2+1)(q-1)/2,(q-3)/2,(q-1)^2/2)$.
The way we construct these graphs in Sage is described in Sect.~\ref{cope05}.
\end{itemize}

\subsection{Graphs from combinatorics}\label{combgraphs}
\begin{itemize}
\item Johnson Graphs $J(m,2)$, see Example~\ref{johnson}.
\item Orthogonal Array block graphs $\OA(k,n)$.
  Sage is able to build a very substantial state-of-the-art collection of orthogonal
  arrays (often abbreviated as OA), thanks to a large implementation project undertaken in
  2013/2014 by the first author in a very productive collaboration with
  Julian R. Abel and Vincent Delecroix.
  For the present work no new constructions of OAs were needed,
  and the link between Sage's OAs and Strongly Regular Graphs
  databases filled in three new entries in Andries Brouwer's database.

\item Steiner Graphs (intersection graphs of BIBD) -- Sage can already build
  several families of Balanced Incomplete Block Designs (when $k\leq 5$, or
  projective planes, or other recursive constructions and fixed-size
  instances). More constructions from \cite{Handbook}
  were added to Sage while working on this project.
\item Goethals-Seidel graphs, see \cite{GS70}.
\item Haemers graphs, see \cite[Sect.8.A]{BvL84}.
\item RSHCD -- graphs from $(n,\epsilon)$-regular symmetric Hadamard matrices $M$ with constant
  diagonal, see Example~\ref{rshcd} for the definition. Several constructions from the
  literature (and one apparently new one, cf. Sect.~\ref{rshcd324})
  for this class of Hadamard matrices were implemented in Sage and are available
  in its Hadamard matrices module.
\item Two-graph descendants. Each \emph{regular two-graph} (a certain class of 3-uniform
$v$-vertex hypergraphs having $2\mu$ three-edges on each pair of points, cf.
e.g. \cite[Chap.10]{BH12}) gives rise to
a \srg\ with parameters $(v-1,2\mu,3\mu-v/2,\mu)$ obtained by descendant construction,
see e.g. \cite[Sect.10.3]{BH12}.
\item Switch $\OA$ Graphs -- these \srg s are obtained from OA
  block graphs (see above). From such a graph $G$ obtained from an $\OA(k,n)$,
  the procedure is to (1) add a new isolated vertex $v$; (2) perform Seidel
  switching on the union of $\{v\}$ and several disjoint $n$-cocliques of
  $G$. Note that a $n$-coclique in $G$ corresponds to a parallel class of the
  $\OA(k,n)$, and that those are easily obtained from an $\OA(k+1,n)$ (i.e. a {\em
    resolvable} $\OA(k,n)$).
\item Polhill Graphs -- In \cite{Polhill09}, Polhill produced 5 new strongly
  regular graphs on 1024 vertices as Cayley graphs. His construction is able to
  produce larger \srg s of order $\geq 4096$, though the current implementation
  only covers the $n=1024$ range.
\item Mathon's pseudo-cyclic \srg s related to symmetric conference matrices,
optionally parameterised by a \srg\ with parameters of a Paley graph, and a skew-symmetric
Latin square \cite{Mat78, ST78}.
\item Pseudo-Paley and Pasechnik graphs from skew-Hadamard matrices.
These are constructions due to Goethals-Seidel \cite{BvL84} and
Pasechnik \cite{Pa92}, constructing graphs on $(4m-1)^2$ vertices
from skew-Hadamard matrices of order $4m$.
Sage builds the corresponding skew-Hadamard matrices from a small
database featuring classical constructions of skew-Hadamard matrices
from \cite{Ha83} and small examples from (anti)-circulant matrices
\cite{GS70s,Wall71}.

\end{itemize}

\subsection{Novel constructions}\label{novel}
Here we collect descriptions of constructions of graphs that are in our view
suffciently novel and interesting to mention. Namely, Sect. \ref{cope05} describes another
construction for a known graph, Sect. \ref{rshcd324} describes a \emph{working} construction
for a graph which was claimed to exist in the literature, although we were unable to verify
a number of published constructions (see Sect. \ref{sect:incorrect} for details).
Finally, Sect. \ref{sect:IK765} discusses an unpublished construction by
Ionin and Kharaghani.

\subsubsection{Cossidente-Penttila hemisystems}\label{cope05}
The construction of the hemisystem in \cite{CP05} requires building $\GQ(q^2,q)$, which is slow.
Thus we designed, following a suggestion of T. Penttila, a more efficient approach,
working directly in $\PG(5,q)$. The partition in question is invariant
under the subgroup $H=\Omega_3(q^2)<O_6^-(q)$. Without loss in generality $H$
leaves the form $B(X,Y,Z)=XY+Z^2$ invariant. We
pick two orbits of $H$ on the $\mathbb{F}_q$-points, one of them $B$-isotropic, with
a representative $(1:0:0)$, viewed as a point of $\Pi:=\PG(2,q^2)$, and
the other corresponding to points of
$\Pi$ that have all the lines on them intersecting the conic of $\Pi$ specified by $B$
in zero or two points. We take $(1:1:\epsilon)$ as a representative, with
$\epsilon\in\mathbb{F}_{q^2}^*$ so that $\epsilon^2+1$ is not a square in $\mathbb{F}_{q^2}$.

Indeed, the conic can be viewed $\{(0:1:0)\}\cup\{(1:-t^2:t)\mid t \in\mathbb{F}_{q^2}\}$.
The coefficients of a generic line on $(1:1:\epsilon)$ are $[1:-1-\epsilon b:b]$,
for $-1\neq \epsilon b$. Thus, to make sure that its intersection with the
conic is always even, we need that the discriminant of
$1+(1+\epsilon b)t^2+tb=0$ never vanishes, and this is if and only if
$\epsilon^2+1$ is not a square.

Finally, we need to adjust $B$, by multiplying it by appropriately
chosen $\nu\in\mathbb{F}_{q^2}^*$, so that $(1:1:\epsilon)$
becomes isotropic under the relative trace norm
$(X:Y:Z)\mapsto \nu B(X,Y,Z)+(\nu B(X,Y,Z))^q$, used to define adjacency in $\Gamma$.

\subsubsection{Regular symmetric Hadamard matrices of order 324.}\label{rshcd324}
We recall the definition of $\RSHCD^+$ and $\RSHCD^-$ from Example~\ref{rshcd}.
An example $M^+$ of $\RSHCD^+$ order 324 was constructed by
Janko, Kharaghani, and Tonchev in \cite{JKT01}, and we implemented
their construction in Sagemath. See \cite{MX06,HX10} for other examples of
$\RSHCD^+$ of order 324.

We use $M^+$ to build an example $M^-$ of $\RSHCD^-$ of order 324,
as follows. One is tempted to apply \cite[Lemma 11]{HX10} to $M^+$,
which says that for an $\RSHCD^\epsilon$ matrix $M$
built from four $n\times n$-blocks $M_{ij}$, so that
\begin{equation}\label{Mpm}
M=\begin{pmatrix} M_{11} & M_{12}\\ M_{21} & M_{22} \end{pmatrix},
\quad\text{the matrix }
T(M):=\begin{pmatrix} M_{11} & -M_{12}\\ -M_{21} & M_{22} \end{pmatrix}
\end{equation}
is an $\RSHCD^{-\epsilon}$, provided that row sums of $M_{11}$ and $M_{22}$ are 0.
However, the latter condition does not hold for $M=M^+$.
We are able to ``twist'' $M^+$ so that the resulting
matrix is amenable to this Lemma. Namely, it turns our that the matrix
\[
    M':=\begin{pmatrix} M_{12} & M_{11}\\ M_{11}^\top & M_{21} \end{pmatrix},
\qquad\text{where } M=M^+,
\]
is $\RSHCD^+$, its diagonal blocks having  row
sums 0, as needed by \eqref{Mpm}. Interestingly, the
$(324,152,70,72)$-\srg\ corresponding to $T(M')$
has a vertex-transitive automorphism group of order 2592, twice the order of the
(intransitive) automorphism group of the $(324,153,72,72)$-\srg corresponding to $M^+$.
As far as we know, this is the only known example of such a vertex-transitive graph.
Other graphs with such parameters were constructed in \cite{HX10}.

\subsubsection{A (765,192,48,48)-graph} \label{sect:IK765}
We were unable to implement the construction of a graph with these parameters
described, as a part of an infinite family, in Ionin and Kharaghani~\cite{IoninKharaghani03}.
The authors of the latter were very kind to send us an updated construction of the
graph in question,
which we successfully implemented in Sagemath, see \cite{trac19990}.
This construction can be found in the documentation of the Sagemath
function {\tt graphs.IoninKharaghani765Graph}.

They have also posted an update \cite{IK16} to \cite{IoninKharaghani03};
we have not yet tried to implement the updated version in full generality.

\section{Missing values}
Among the 1150 realizable, according to Andries Brouwer's
database, parameter sets,  our implementation can realize 1142. Up to taking graph complements,
the list of currently missing entries is as follows.
\begin{center}
\begin{tabular}{ccccl}
(196&90&40&42)& $\RSHCD^-$~(may not exist, cf. Sect.~\ref{sect:incorrect})\\
(196&135&94&90)& Huang, Huang and Lin \cite{HHL09}\\
(378&116&34&36)& Muzychuk $S6(n=3,d=3)$~\cite{Muzychuk07}\\
(512&133&24&38)& Godsil~$(q=8,r=3)$~\cite{Godsil92}\\
\end{tabular}
\end{center}
The fisrt entry is discussed in Sect.~\ref{sect:incorrect}. Implementation of the remaining three
entries is currently in progress.

\section{Incorrect RSHCD constructions ?}\label{sect:incorrect}
We were unable to reproduce the following two constructions of Regular Symmetric Hada\-mard matrices
with Constant Diagonal (RSHCDs) and thus the corresponding \srg s (in the sense of
\cite[Sect.10.5]{BH12}).
\begin{itemize}
\item In \cite[Sect.10.5.1, (i)]{BH12}, the construction of $\RSHCD(196,-)$ is attributed to
  \cite{IS06}, in which the existence of a $(4k^2,2k^2+k,k^2+k)$-strongly
  regular graph, equivalent to a $\RSHCD(196,-)$ for $k=7$, is claimed in Theorem
  8.2.26.(iii). The latter says that the $\RSHCD(196,-)$ can
  be easily obtained from the $\RSHCD(196,+)$ from \cite[Theorem
  8.2.26.(ii)]{IS06}. While the construction of (ii) was successfully
  implemented in Sage, following the authors' instructions for (iii) did not
  lead us to the $\RSHCD(196,-)$. Communication with the authors did not solve
  the issue, and we are not aware of any other proof of the existence of a
  $(196,90,40,42)$-\srg.

\item In \cite[Sect.10.5.1, (iii)]{BH12} one finds the following claim,
  attributed to \cite[Corollary 5.12]{SWW72}.
  \begin{center}
    {\it If $n-1$ and $n+1$ are odd prime powers, there exists a $\RSHCD(n^2,+)$.}
  \end{center}
  We implemented the construction provided in~\cite[Corollary 5.12]{SWW72}, but
  that did not lead us to the expected strongly regular graph. We also note that
  while Corollary 5.12 does not claim that the provided matrices are regular, that
  claim appears in the theorem on which it relies. The author of~\cite[Corollary
  5.12]{SWW72} did not answer our message, and we discarded this construction as
  broken in our work.

  This construction should have been able to produce a $\RSHCD(676,+)$ and a
  $\RSHCD(900,+)$. Fortunately in the end it was not required.

\end{itemize}

The following construction of RSHCDs needed a lot of effort and a number of discussions with
Andries Brouwer to correct crucial misprints in several sources and combine them into a working
construction.
\begin{itemize}
\item In \cite[Sect.10.5.1, (iv)]{BH12} one finds the following claim,
  attributed to \cite[Corollary 5.16]{SWW72}.

  \begin{center}
    {\it If $a+1$ is a prime power and there exists a symmetric conference matrix of order $a$, then there exists a $\RSHCD(a^2,+)$.}
  \end{center}
  Following [loc.cit.] did not lead us to the expected result; as it
  turns out, [loc.cit.] has a typo, and the correct formulae should be taken from
  the original source \cite[Corollary~17]{MR0316273} by J. Wallis and A.\,L. Whiteman.
  An essential ingredient in this construction, referred to in \cite{MR0316273},
  is a special pair of difference sets due to G. Szekeres \cite[Theorem~16]{MR0246787},
  defined in (4.1) and (4.2) there. However, (4.2) has a typo (it has $-$ instead of $+$ sign),
  invalidating the construction. Fortunately, a correct definition for may be found
  \cite[Theorem 2.6, p.303]{SWW72}. This construction allowed us to produce
  $\RSHCD(676,+)$ and a $\RSHCD(900,+)$.
\end{itemize}

\section{Entries added to database during this work}

By linking Sage's database of Orthogonal Arrays with its database of
Strongly Regular Graphs, we were able to fill in the following three values:
\begin{itemize}
\item $(196, 78, 32, 30)$ -- can be obtained from an $\OA(6,14)$~\cite{Todorov12}\\[-6mm]
\item $(324, 102, 36, 30)$ -- can be obtained from an $\OA(6,18)$~\cite{Abel15}\\[-6mm]
\item $(324, 119, 46, 42)$ -- can be obtained from an $\OA(7,18)$~\cite{Abel15}
\end{itemize}

This can be seen as a by-product of making two mathematical databases,
which formerly only existed in printed form, inter-operable.
In our implementation, any update of the combinatorial designs databases can be beneficial for
the database of strongly regular graphs.\\

We obtained a $(1024,462,206,210)$-graph while going through
the constructions from~\cite{Polhill09}, although this value did not appear in the
online database at that time.

\section*{Acknowledgments}
The authors thank Andries Brouwer, Eric Chen, Luis Disset, Hadi Kharaghani,
Misha Muzychuk, Tim Penttila, John B. Polhill, Leonard Soicher, Vladimir Tonchev,
and Alfred Wasserman for many
helpful discussions and communications.

The second author was partially  supported by the EU Horizon 2020 research
and innovation programme, grant agreement OpenDreamKit No 676541.

\bibliography{db}
\bibliographystyle{abbrv}
\end{document}